\documentstyle[12pt,twoside,leqno]{amsart}

\newtheorem{Theorem}{Theorem}[section]

\newtheorem{Lemma}[Theorem]{Lemma}

\newtheorem{Question}[Theorem]{Question}

\theoremstyle{remark}
\newtheorem{Remark}[Theorem]{Remark}

\theoremstyle{definition}
\newtheorem{Definition}[Theorem]{Definition}

\begin{document}
\bibliographystyle{plain}
\title{Strong $A_{\infty}$-weights, Besov and Sobolev capacities in metric measure spaces}
\author{\c Serban Costea}
\thanks{Work partially supported by NSERC, by the Emil Aaltonen foundation, by the Fields Institute,
 and by the NSF grant DMS 0244421}

\address{McMaster University, Department of Mathematics and Statistics, 1280 Main Street West,
Hamilton, Ontario L8S 4K1, CANADA}

\address{Fields Institute for Research in Mathematical Sciences, 222 College Street, Toronto, Ontario M5T 3J1, CANADA}
\email{secostea@@math.mcmaster.ca}
\begin{abstract}

 This article studies strong $A_{\infty}$-weights in Ahlfors $Q$-regular and geodesic metric spaces satisfying a weak
 $(1,s)$-Poincar\'{e} inequality for some $1<s \le Q<\infty.$ It is shown that whenever $\max(1, Q-1)<s \le Q,$ a function
 $u$ yields a strong $A_{\infty}$-weight of the form $w=e^{Qu}$ if the minimal $s$-weak upper gradient $g_{u}$ has sufficiently
 small $||\cdot||_{{\mathcal{L}}^{s, Q-s}(X,\mu)}$ norm. Similarly, it is proved that if $1<Q<p<\infty,$ then $w=e^{Qu}$ is a strong $A_{\infty}$-weight whenever the Besov $p$-seminorm $[u]_{B_p(X)}$ of $u$ is sufficiently small.
\end{abstract}

\maketitle
\section{Introduction}
 In this paper $(X,d,\mu)$ is a complete and unbounded metric measure space. In addition, we assume that it is Ahlfors $Q$-regular for some
 $Q>1.$ That is, there exists a constant $C=c_{\mu}$ such that, for
 each $x \in X$ and all $r>0,$
 $$C^{-1} r^Q \le \mu(B(x,r)) \le C r^Q.$$
 Furthermore, $X$ is assumed to be geodesic. That is, every pair of points can be joined by a curve whose length is the distance between the points.

 We will also assume that $(X,d,\mu)$ satisfies a weak $(1,s)$-Poincar\'{e} inequality for some $s \in (1,Q].$ That is, there exist constants
 $C>0$ and $\lambda \ge 1$ such that for all balls $B$ with radius $r,$ all measurable functions $u$ on $X$ and all upper gradients $g$ of $u$
 we have
 \begin{equation} \label{defn of 1s PI}
 \frac{1}{\mu(B)} \int_{B} |u-u_{B}| \, d\mu \le Cr \left( \frac{1}{\mu(\lambda B)} \int_{\lambda B} g^{s} \, d\mu \right)^{1/s},
 \end{equation}
 where $\lambda B$ represents the ball concentric with $B$ with radius $\lambda$ times the radius of $B$ whenever
 $\lambda>0,$ and $u_{E}$ denotes the average of $u$ on the measurable set $E \subset X$ with respect to the measure $\mu$ whenever $0<\mu(E)<\infty.$
 We recall that a nonnegative Borel function $g$ is an \textit{upper gradient} for a real-valued measurable function
 $u$ on $X$ if for all rectifiable paths $\gamma: [0, l_{\gamma}] \rightarrow X$ we have
 \begin{equation} \label{defn of upper gradient}
 |u(\gamma(0))-u(\gamma(l_{\gamma}))| \le \int_{\gamma} g \, ds.
 \end{equation}
 Here and throughout the paper the rectifiable curve $\gamma:[0, l_{\gamma}] \rightarrow X$ is assumed to be parametrized by the
 arc length $ds,$ where $l_{\gamma}$ is the length of $\gamma.$

 We study sufficient conditions under which we get strong $A_{\infty}$-weights in $X.$ A nontrivial doubling measure $\nu$ on $X$ is a
 Radon measure for which there exists a constant $C>1$ such that
 $$0< \nu(2B) \le C \nu(B)$$
 for all balls $B.$

 To every doubling measure $\nu$ on $X$ we can associate a quasidistance on $X$ defined by
 \begin{equation}\label{def delta nu}
 \delta_{\nu}(x,y)=\nu(B_{xy})^{1/Q},
 \end{equation}
 where $B_{xy}:=B(x, d(x,y)) \cup B(y, d(y,x)).$
 To say that $\delta_{\nu}(x,y)$ is a quasidistance means by definition that $\delta_{\nu}: X \times X \rightarrow [0, \infty)$ is
 symmetric, vanishes if and only if $x=y,$ and satisfies
 \begin{eqnarray} \label{def quasidistance}
 \delta_{\nu}(x,z) \le C(\delta_{\nu}(x,y)+\delta_{\nu}(y,z))
 \end{eqnarray}
 for some $C \ge 1$ and all $x,y,z \in X.$ If (\ref{def quasidistance}) was satisfied with $C=1,$ then the quasidistance $\delta_{\nu}$
 would in fact be a distance function.

  We call $\nu$ a \textit{metric doubling measure} if the quasidistance $\delta_{\nu}$
 is comparable to a distance $\delta'_{\nu};$ that is, there exists a
 distance function $\delta'_{\nu}$ on $X$ and a constant $C>0$
 such that
 \begin{equation}\label{quasidistance comparable to a distance}
 C^{-1} \delta_{\nu}(x,y) \le \delta'_{\nu}(x,y) \le C \delta_{\nu}(x,y) \mbox{ for all } x,y \in X.
 \end{equation}

 We say that a nonnegative function $w \in L_{loc}^{1}(X)$ is an \textit{$A_{p}$-weight} with respect to the measure $\mu$ for some $1 < p <\infty$ and we write $w \in A_{p}(\mu)$ if there exists a constant $C \ge 1$ such that
 \begin{equation*}
 \left(\frac{1}{\mu(B)} \int_{B} w(x)^{-1/(p-1)} \, d\mu(x)\right)^{p-1} \frac{1}{\mu(B)} \int_{B} w(x) \, d\mu(x) \le C
 \end{equation*}
 for all balls $B \subset X.$
 We say that $w$ is a $A_{\infty}$-weight with respect to the measure $\mu$ and we write $w \in A_{\infty}(\mu)$ if $w$ is an $A_p$-weight with respect to $\mu$ for some $p$ in $(1,\infty).$ That is, $$A_{\infty}(\mu) = \cup_{p>1} A_p(\mu).$$

 We define $w$ to be a strong $A_{\infty}$-weight if it is the density of a metric doubling measure $\nu$ and moreover,
 it is an $A_{\infty}$-weight with respect to $\mu.$ That is,
 $$d\nu(x)=w(x) \, d\mu(x)$$
 where $w \in A_{\infty}(\mu)$ and $\nu$ is a metric doubling measure.

 Strong $A_{\infty}$-weights in ${\mathbf{R}}^n$ were introduced in the early 90's by David and Semmes
 in \cite{DS} and \cite{Sem} when trying to identify the subclass of $A_{\infty}$-weights that are comparable
 to the Jacobian determinants of quasiconformal mappings.
 \begin{Question} In the Euclidean setting, metric doubling measures have densities that are $A_{\infty}$-weights.
 (See \cite{Sem}.) An open question in the metric setting is whether or not metric doubling measures
 necessarily have $A_{\infty}$-densities.
 \end{Question}
 In the last few years strong $A_{\infty}$-weights were studied by Bonk, Heinonen, and Saksman in \cite{BHS1}
 and \cite{BHS2} and by the author in \cite{Cos1}.

 In the Euclidean setting Bonk and Lang proved in \cite{BL} that if $\nu$ is a signed Radon measure on ${\mathbf{R}}^2$
 such that $\nu^{+}({\mathbf{R}}^2)<2 \pi$ and
 $\nu^{-}({\mathbf{R}}^2)<\infty,$ then $({\mathbf{R}}^2, \widetilde D_{\nu})$
 is bi-Lipschitz equivalent to ${\mathbf{R}}^2$ endowed with the
 Euclidean metric, where
 \begin{equation*}
 \widetilde D_{\nu}(x,y)=\inf \bigg \lbrace \int_{\alpha} e^{u} ds:
 \alpha \mbox { analytic curve connecting $x,y$}  \bigg \rbrace,
 \end{equation*}
 $u$ is a solution of $-\Delta u=\nu$ with $|\nabla u| \in L^{2}({\mathbf{R}}^2),$ and
 $\nu=\nu^{+}-\nu^{-}$ is the Jordan decomposition of $\nu.$ In particular, it is
 proved that $w=e^{2u}$ is comparable to the Jacobian of a quasiconformal mapping
 $f: {\mathbf{R}}^2 \rightarrow {\mathbf{R}}^2,$ which implies that $w$ is a strong
 $A_{\infty}$-weight.

 Here we prove a result in $(X, d, \mu),$ related to \cite[Theorem 5.1]{Cos1} and to the result from \cite{BL}.
 It states that $A_{\infty}$-weights of the form $w=e^{Qu}$ are strong $A_{\infty}$-weights if $u$
 is a locally integrable function that has an upper gradient $g$ in the Morrey space ${{\mathcal L}}^{s,Q-s}(X, \mu)$
 with small $||\cdot||_{{{\mathcal L}}^{s,Q-s}(X, \mu)}$ norm for some $s>1$ lying in $(Q-1,Q].$

 We say that for $1 \le s \le Q,$ the Morrey space ${{\mathcal L}}^{s,Q-s}(X, \mu)$ is defined to be the linear space
 of locally $\mu$-integrable functions $u$ on $X$ such that
 $$||u||_{{{\mathcal L}}^{s,Q-s}(X, \mu)} = \sup_{x \in X} \sup_{r>0} \left( r^{s-Q} \int_{B(x,r)} |u(y)|^{s} \, d\mu(y) \right)^{1/s}.$$
 In particular ${{\mathcal L}}^{Q,0}(X, \mu)=L^Q(X).$ We refer to \cite[p.\ 65]{Gia} for more information about Morrey spaces in the Euclidean setting and their use in the theory of partial differential equations.

 If $(X,d,\mu)$ is an Ahlfors $Q$-regular metric space with $Q>1$ satisfying a weak $(1,s)$-Poincar\'{e} inequality for some
 $s \in (1,Q],$ it follows from (\ref{defn of 1s PI}) that there exists a constant $C$ depending only $s$ and on data of $X$ such that
 \begin{equation}\label{BMO seminorm of u dominated by Morrey norm of grad u}
 [u]_{{\mathrm{BMO}}(X)} \le C ||g||_{{{\mathcal L}^{s, Q-s}}(X, \mu)}
 \end{equation}
 whenever $g$ is an upper gradient of $u.$ Here and throughout the paper $[u]_{{\mathrm{BMO}}(X)}$ is the \textit{bounded mean
 oscillation} seminorm that measures the oscillation of $u$ on balls
 in $X,$ given by
 $$[u]_{{\mathrm{BMO}}(X)}=\sup_{a \in X} \sup_{r>0}
 \frac{1}{\mu(B(a,r))} \int_{B(a,r)} |u(x)-u_{B(a,r)}| \, d\mu(x).$$

 In \cite[Theorem 3.1]{BHS1} the authors prove that if $u$ belongs to the
 Bessel potential space $L^{\alpha, \frac{n}{\alpha}}({\mathbf{R}}^n),$ $0<\alpha<n,$
 then $w=e^{nu}$ is a strong $A_{\infty}$-weight with data depending only on $\alpha,$
 $n,$ and the $L^{\alpha, \frac{n}{\alpha}}$-norm of $u.$ Here we prove a result similar to
 \cite[Theorem 3.1]{BHS1} and \cite[Theorem 5.2]{Cos1}. This result yields strong $A_{\infty}$-weights of the form
 $w=e^{Qu}$ when $u$ has small Besov $p$-seminorm, $1<Q<p<\infty.$

 For $1<Q<p<\infty$ we define
 \begin{equation} \label{Besov space defn}
 B_p(X)=\{u \in L^p(X): ||u||_{B_p(X)} <\infty\},
 \end{equation}
 where
 \begin{equation} \label{Besov norm defn}
 ||u||_{B_p(X)}=||u||_{L^p(X)}+[u]_{B_p(X)}
 \end{equation}
 with
 \begin{equation} \label{Besov seminorm defn}
 [u]_{B_p(X)}=\left( \int_{X} \int_{X}
 \frac{|u(x)-u(y)|^{p}}{d(x,y)^{2Q}} \, d\mu(x) \, d\mu(y) \right)^{1/p}.
 \end{equation}

 The expressions $||u||_{B_p(X)}$ and $[u]_{B_p(X)}$ from (\ref{Besov norm defn}) and (\ref{Besov seminorm defn}) are
 called the \textit{Besov $p$-norm} and the \textit{Besov $p$-seminorm} of $u$ respectively.
 If $(X,d,\mu)$ is Ahlfors $Q$-regular, there exists a constant $C$ depending on $p$ and on the data of $X$ such that
 \begin{equation} \label{BMO seminorm of u dominated by Besov seminorm of u}
 [u]_{{\mathrm{BMO}}(X)} \le C [u]_{B_p(X)}
 \end{equation}
 whenever $u \in L_{loc}^{1}(X).$

 Besov spaces have been studied in the last decades by Jonsson and Wallin in \cite{JW}, by Fukushima and Uemura in
 \cite{FU}, by Xiao in \cite{Xia}, and by the author in \cite{Cos1} and \cite{Cos2}.
 Recently they have been used in the study of quasiconformal mappings
 in metric spaces and in geometric group theory. See \cite{Bou} and \cite{BP}.

 Capacities associated with Besov spaces were studied by Netrusov in \cite{Net1} and
 \cite{Net2}, by Adams and Hurri-Syrj\"{a}nen in \cite{AHS}, by Adams and Xiao in
 \cite{AX1} and \cite{AX2}, and by the author in \cite{Cos1}.
 Bourdon in \cite{Bou} and the author in \cite{Cos2} studied Besov
 $p$-capacity in metric settings.

 {\bf{Acknowledgements.}} This article was written when the author was a visiting member of the Fields Institute in 2008.
 Part of the research was done when the author was a visiting postdoctoral researcher at the Helsinki University of Technology in 2007.
 The author wishes to thank Mario Bonk for useful conversations leading to Lemma \ref{Proj lemma} and
 Carlos P\'erez for helpful discussions regarding the John-Nirenberg lemma.

\section{Preliminaries}

 In this section we recall standard definitions and results. The open ball with center $x \in X$ and radius $r>0$ is denoted
 $B(x,r)=\{ y \in X: d(x,y)<r \},$ the closed ball by $\overline{B}(x,r)=\{ y \in X: d(x,y) \le r \},$ and the sphere
 by $S(x,r)=\{ y \in X: d(x,y)=r \}.$
 Throughout this paper, $C$ will denote a positive constant whose value
 is not necessarily the same at each occurrence; it may vary even within
 a line. $C(a,b, \ldots)$ is a constant that depends only on the parameters
 $a,b, \ldots.$
 Here $\Omega$ will denote a nonempty open subset of $X.$ For $E \subset X,$
 the closure and the complement of $E$ with respect
 to $X$ will be denoted by $\overline{E}$ and $X \setminus E$
 respectively; $\mbox{diam }E$ is the diameter of $E$ with
 respect to the metric $d$ and $E \subset \subset F$ means that $\overline{E}$
 is a compact subset of $F.$

 For a measurable $u: \Omega \rightarrow \mathbf{R},$ $\mbox {supp } u$ is the
 smallest closed set such that $u$ vanishes on the complement of $\mbox {supp } u.$
 We also use the spaces
 \begin{eqnarray*}
 Lip(\Omega)&=&\{\varphi : \Omega \rightarrow {\mathbf{R}} : \varphi \mbox{ is Lipschitz} \},\\
 Lip_{0}(\Omega)&=&\{\varphi : \Omega \rightarrow {\mathbf{R}} : \varphi \mbox{ is Lipschitz and supp $\varphi \subset \subset \Omega$} \}.
 \end{eqnarray*}

\subsection{Newtonian spaces}

We introduce now some definitions and known results about Newtonian spaces to be used in this paper.
Let $1\le s<\infty.$ The $s$-modulus of a family of paths $\Gamma$ in $X$ is the number
$$\inf_{\rho} \int_{X} \rho^s \, d\mu,$$
where the infimum is taken over all non-negative Borel measurable functions $\rho$ such that for all
rectifiable paths $\gamma$ which belong to $\Gamma$ we have
$$\int_{\gamma} \rho \, ds \ge 1.$$
It is known that the $s$-modulus is an outer measure on the collection of all paths in $X.$

A property is said to hold for $s$-\textit{almost all} paths, if the set of paths for which the property fails
is of zero $s$-modulus. If (\ref{defn of upper gradient}) holds for $s$-almost all paths $\gamma,$ then $g$ is said to be
a $s$-weak upper gradient of $u.$ We could have stated the definition of the weak $(1,s)$-Poincar\'{e} inequality
by requiring the inequality (\ref{defn of 1s PI}) to hold for all $s$-weak upper gradients of $u.$ (See \cite{KoM}.)
Similarly we can define weak $(q,s)$-Poincar\'{e} inequalities for $q>1.$

\textit{Without further notice, we assume that $1 < s <\infty.$} We define the space $\widetilde{N}^{1,s}(X)$ to be the collection of all the functions $u$ that are
$s$-integrable and have a $s$-integrable $s$-weak upper gradient $g.$ This space is equipped with the norm
$$||u||_{\widetilde{N}^{1,s}(X)}=\left(||u||_{L^s(X)}^s+ \inf ||g||_{L^s(X)}^s\right)^{1/s},$$
where the infimum is taken over all $s$-weak upper gradients of $u.$ The \textit{Newtonian space} on $X$ is the quotient space
 $$N^{1,s}(X)=\widetilde{N}^{1,s}(X)/\sim$$
with the norm $||u||_{N^{1,s}(X)}=||u||_{\widetilde{N}^{1,s}(X)},$ where $u \sim v$ if and only if $||u-v||_{\widetilde{N}^{1,s}(X)}=0.$
For basic properties of the Newtonian spaces we refer to \cite{Sha1}. Cheeger in \cite{Che} gives an alternative definition
which leads to the same space when $1<s<\infty.$ For future reference we recall some known facts (see \cite{KiM} and \cite{Sha2}):

(i) The functions in $\widetilde{N}^{1,s}(X)$ are defined outside a path family of $s$-modulus zero. This implies that the functions
in $\widetilde{N}^{1,s}(X)$ cannot be changed arbitrarily on sets of measure zero.

(ii) If $1<s<\infty,$ every function $u$ that has a $s$-integrable $s$-weak upper gradient has in fact a minimal $s$-integrable $s$-weak upper gradient in $X,$ denoted by $g_{u},$ in the sense that if $g$ is another $s$-weak upper gradient of $u,$ then $g_{u} \le g$ $\mu$-a.e.\
in $X.$

(iii) For every $c \in {\mathbf{R}}$ the minimal $s$-weak upper gradient satisfies $g_{u}=0$ $\mu$-a.e.\ on the set $\{ x \in X: u(x)=c \}.$

(iv) If $u \in N^{1,s}(X)$ and $v$ is a bounded Lipschitz continuous function, then $uv \in N^{1,s}(X)$ and $g_{uv} \le |u|g_{v}+|v|g_{u}$
$\mu$-a.e.

We emphasize that these properties hold without any additional assumptions on the measure $\mu$ and on the space $X.$

The $s$-capacity of a set $E \subset X$ is defined by (see \cite{BBS})
$$C_s(E)= \inf_{u} ||u||_{N^{1,s}(X)}^s,$$
where the infimum is taken over all functions $u \in N^{1,s}(X)$ whose restriction on $E$ is bounded below by $1.$
A property is said to hold $s$-\textit{quasieverywhere} (or $s$-q.e.), if it holds everywhere except on a set of $s$-capacity zero. A function is
$s$-\textit{quasicontinuous}, if there is an open set of arbitrarily small $s$-capacity such that the function is continuous when
restricted to the complement of the set. Every function in $\widetilde{N}^{1,s}(X)$ is defined $s$-quasieverywhere. Moreover, if
$u,v \in N^{1,s}(X)$ and $u=v$ $\mu$-a.e.,\ then $u=v$ $s$-quasieverywhere. In particular, this implies that $u$ and $v$ belong to
the same equivalence class in $N^{1,s}(X).$

We introduce the notion of a local Newtonian space as follows.
\begin{Definition} \label{defn for local newtonian spaces}
 We say that $u$ belongs to the \textit{local Newtonian space} $N_{loc}^{1,s}(X)$
if $u \in N^{1,s}(\Omega)$ for every open set $\Omega \subset \subset X.$ If $u \in N_{loc}^{1,s}(X)$ with $1<s<\infty,$ then $u$ has a minimal
$s$-weak upper gradient $g_{u}$ in $X$ in the following sense: if $\Omega \subset \subset X$ is an open set and $g$ is the minimal upper gradient
of $u$ in $\Omega,$ then $g_{u}=g$ $\mu$-a.e.\ in $\Omega.$
\end{Definition}

 From now on throughout the rest of the paper we assume that the measure $\mu$ is Borel and Ahlfors $Q$-regular for some $Q>1.$ Furthermore we assume that the space supports a weak $(1,s)$-Poincar\'{e} inequality for some $1<s \le Q.$ We recall a few useful properties of Newtonian spaces that hold under these additional assumptions (see \cite{BBS} and \cite{KiM}):

 (i) The space $X$ is proper (that is, closed and bounded sets are compact).

 (ii) Lipschitz functions are dense in $N^{1,s}(X)$ and Lipschitz functions which vanish in the complement of an open set $\Omega$ are dense in $N_{0}^{1,s}(\Omega),$ where
 $$N_{0}^{1,s}(\Omega)=\{ u \in N^{1,s}(X): u=0 \mbox{ $s$-q.e.\ in $X \setminus \Omega$}  \}.$$

 (iii) Every function in $N^{1,s}(X)$ is $s$-quasicontinuous.

Now we introduce the relative Sobolev $s$-capacity as in \cite{Cos3}. See also \cite{Bjo}.

\begin{Definition} \label{defn for relative sobolev s capacity}
Let $1<s,Q<\infty.$ Suppose $(X,d,\mu)$ is a proper and unbounded Ahlfors $Q$-regular metric space that
satisfies a weak $(1,s)$-Poincar\'{e} inequality. Let $\Omega \subset \subset X$ be open. For $E \subset \Omega$ we let
$$A(E, \Omega) = \{ u \in N_{0}^{1,s}(\Omega): u \ge 1 \mbox{ in a neighborhood of $E$} \}.$$ We call $A(E, \Omega)$ the
set of \textit{admissible functions for the condenser} $(E, \Omega).$ The relative $s$-capacity of the pair $(E, \Omega)$
is defined by
$${\mathrm{cap}}_{s}(E, \Omega) = \inf \bigg \lbrace \int_{\Omega} g_{u}^s \, d\mu: u \in A(E, \Omega)   \bigg \rbrace.$$
\end{Definition}

\subsection{Besov spaces and capacities} Now we introduce some definitions and results about Besov spaces and capacities
to be used in this paper. We follow \cite{Cos2}. See also \cite{Cos1}.

Let $1<Q<p<\infty$ be fixed. Suppose $(X,d,\mu)$ is an Ahlfors $Q$-regular metric space. For an open set $\Omega \subset X$ we define
$$B_p(\Omega)= \{ u \in B_p(X): u=0 \mbox{ $\mu$-a.e.\ in  $X \setminus \Omega$}  \},$$
where $B_p(X)$ is defined as in (\ref{Besov space defn}). For a function $u \in B_p(\Omega)$ we let
$$||u||_{B_p(\Omega)}=||u||_{B_p(X)} \mbox{ and } [u]_{B_p(\Omega)}=[u]_{B_p(X)}.$$
We notice that $Lip_{0}(\Omega) \subset B_p(\Omega)$ when $1<Q<p<\infty.$ We define $B_p^{0}(\Omega)$ as the closure of $Lip_{0}(\Omega)$
in $B_p(\Omega)$ with respect to the Besov $p$-norm. It has been proved in \cite{Cos2} that $B_p(X),$ $B_p(\Omega),$ and $B_p^{0}(\Omega)$
are reflexive spaces. (See \cite[Lemma 3.1]{Cos2} and the discussion before \cite[Lemma 3.4]{Cos2}.)

The Besov $p$-capacity of a set $E \subset X$ is defined by (see \cite{Cos2})
$${\mathrm{Cap}}_{B_p}(E) = \inf \{ ||u||_{L^p(X)}^p+ [u]_{B_p(X)}^p \},$$
where the infimum is taken over all functions $u \in B_p(X)$ that are bounded from below by $1$ in an open neighborhood of $E.$
A property is said to hold Besov $p$-\textit{quasieverywhere} (or simply $B_p$-q.e.), if it holds everywhere except a set of Besov
$p$-capacity zero. A locally integrable function $u$ is called \textit{$B_p$-quasicontinuous} if there exists an
open set of arbitrarily small Besov $p$-capacity such that $u$ is continuous when restricted to the complement of the set.

\begin{Remark} \label{remark about Besov quasicontinuity of representatives}
 It has been shown in \cite{Cos2} that if $u \in B_p(X),$ then there exists a $B_p$-quasicontinuous function $v$ such that
 $u=v$ $\mu$-a.e. Such a function $v$ is called a \textit{quasicontinuous representative} of $u.$ In addition, we can choose
 $v$ to be Borel. Moreover, two such quasicontinuous representatives agree in fact $B_p$-q.e. Similar statements were proved
 if $u \in L_{loc}^{1}(X)$ with $[u]_{B_p(X)}<\infty.$ (See \cite[Section 5]{Cos2}.)
\end{Remark}

Suppose $\Omega \subset X$ is open. For $E \subset \Omega$ the relative Besov $p$-capacity of the condenser $(E, \Omega)$ is defined by
(see \cite{Cos2})
$${\mathrm{cap}}_{B_p}(E, \Omega) = \inf \{ [u]_{B_p(\Omega)}^p: u \in B_p^{0}(\Omega) \mbox{ and $u \ge 1$ in a neighborhood of $E$} \}.$$

\section{Main results}

In this section we present the results about strong $A_{\infty}$-weights. We prove the following theorems.

\begin{Theorem} \label{SA-infty weight s<Q}
Let $1<s \le Q<\infty$ be fixed. We assume that $s>Q-1.$ Suppose $(X,d,\mu)$ is an Ahlfors $Q$-regular and
geodesic unbounded metric space satisfying a weak $(1,s)$-Poincar\'{e} inequality. Let $u \in N_{loc}^{1,s}(X)$
be such that it has a minimal $s$-weak upper gradient $g_{u}$ in the Morrey space ${\mathcal{L}}^{s, Q-s}(X,\mu).$
There exists a constant $\varepsilon>0$ depending only on $s$ and on the data of $X$ such that if
$$||g_{u}||_{{\mathcal{L}}^{s, Q-s}(X,\mu)}<\varepsilon,$$
then $w=e^{Qu}$ is a strong $A_{\infty}$-weight with data depending only on $s$ and on the data associated with $X.$

\end{Theorem}

\begin{Theorem} \label{SA-infty weight p>Q} Let $1<s<Q<p<\infty$ be fixed. Suppose $(X,d,\mu)$ is an Ahlfors $Q$-regular and geodesic
unbounded metric space satisfying a weak $(1,s)$-Poincar\'{e} inequality. Let $u \in L_{loc}^{1}(X)$
be such that $[u]_{B_p(X)}<\infty.$ There exists a constant $\varepsilon>0$ depending only on $p$ and on the data of $X$ such that if
$$[u]_{B_p(X)}<\varepsilon,$$
then $w=e^{Qu}$ is a strong $A_{\infty}$-weight with data depending only on $p$ and on the data associated with $X.$

\end{Theorem}

For $r \in (0, \infty)$ we define the Hausdorff $r$-content of a set $E \subset X$ by
$$\Lambda_r^{\infty}(E) = \inf \{ \sum_{i} \mbox{ diam}(G_i)^r: E \subset \bigcup_i G_i \},$$
where the infimum is taken over all coverings of $E$ by open sets $G_i.$

The following lemma is a generalization of \cite[Lemma 3.11]{BHS1}. We again thank Mario Bonk
for his contribution to this result.

\begin{Lemma} \label{Proj lemma} Suppose $(X,d,\mu)$ is a
proper and unbounded geodesic Ahlfors $Q$-regular metric space
admitting a weak $(1,Q)$-Poincar\'{e} inequality for some
$1<Q<\infty.$ Suppose $0<r\le 1.$ Let $x,y \in X$ and let $E
\subset X$ be a bounded Borel set. Suppose that $B_{1}, \ldots,
B_{k}$ are open balls such that $x \in B_{1}, y \in B_{k}$ and
$B_{i} \cap B_{i+1} \neq \emptyset$ for $i=1,\ldots, k-1.$ Then
there exists constants $c_1$ and $C$ depending on $r$ and
on the data of $X$ with the following property: if
\begin{equation} \label{Proj Lemma Lambdar le c1 distr}
\Lambda_{r}^{\infty}(E) \le c_1 \, d(x,y)^r,
\end{equation}
then
\begin{equation} \label{Proj Lemma sum dominates distr}
\sum_{i \in {\mathcal{G}}_0} {\rm{diam}}(B_i)^r > \frac{1}{(20C)^r} d(x,y)^r,
\end{equation}
where
\begin{equation} \label{Proj Lemma defn of G0}
{\mathcal{G}}_0= \bigg \lbrace i=1, \ldots, k: \mu(E \cap
B_i) \le \frac{1}{2} \mu(B_i) \bigg \rbrace.
\end{equation}
\end{Lemma}

\begin{pf} Since $X$ is Ahlfors $Q$-regular, proper and geodesic, it
follows that it is also locally linearly connected. That is, there exists
a constant $C \ge 2$ such that every pair of points in $B(x,R)$ can be joined
by a rectifiable path in $B(x,CR)$ and every pair of points in $X \setminus B(x,R)$ can be joined
by a continuum in $X \setminus B(x, R/C).$ (See \cite[Section 3]{HeK} and
\cite[Sections 8,9]{Hei}.)

We choose a family ${\mathcal{I}} \subset \{1, \ldots,k  \}$ such that
$$CB_i \cap CB_j = \emptyset \mbox{ whenever } i \neq j \in {\mathcal{I}} \mbox { and } \bigcup_{i=1}^{k}
CB_i \subset \bigcup_{i \in {\mathcal{I}}} 5CB_i,$$
where $C$ is the constant associated with the locally linear connectivity of $X.$
(See \cite[Theorem 1.16]{Hei}.) For every $i=1,\ldots,k-1$
let $x_i \in B_{i} \cap B_{i+1}$. We let $x_0=x$ and $x_k=y.$ Since
$X$ is locally linearly connected, we have that for every $i=1, \ldots, k$ there exists a
rectifiable path $\gamma_i$ in $CB_i$ connecting $x_{i-1}$ and
$x_i.$ This yields a rectifiable path $\gamma$ in $\cup_{i=1}^{k}
CB_{i}$ connecting $x$ and $y$ and therefore
\begin{equation} \label{upper bound for Lambdas of I}
\Lambda_r^{\infty}(\bigcup_{i \in {\mathcal{I}}} 5CB_{i}) \ge
\Lambda_r^{\infty}(\bigcup_{i=1}^{k} CB_{i}) \ge d(x,y)^r.
\end{equation}

We can assume without loss of generality that
$\Lambda_{r}^{\infty}(E)>0.$ Let $(D_j)_{j \in {\mathcal{J}}}$ be a
countable covering by open balls for $E$ such that
$$\frac{1}{5} D_i \cap \frac{1}{5} D_j = \emptyset $$
whenever $i \neq j \in {\mathcal{J}}$ (see \cite[Theorem 1.16]{Hei}) and such that
$$\sum_{j \in {\mathcal{J}}} \mbox{diam}(D_j)^r < 2^{r+1} \Lambda_{r}^{\infty}(E).$$
For every $i \in {\mathcal{I}}$ we define $${\mathcal{F}}_i = \{ j
\in {\mathcal{J}}: D_j \cap CB_i \neq \emptyset  \}.$$ We denote
$${\mathcal{G}}= \{i \in {\mathcal{I}}: {\rm{diam}}(D_j) \le {\rm{diam}}(CB_i) \, \mbox {for all } j \in
{\mathcal{F}}_i \} \mbox{ and } {\mathcal{B}}={\mathcal{I}} \setminus {\mathcal{G}}.$$
Suppose $i \in {\mathcal{B}}.$ Then there
exists $j=j_i \in {\mathcal{F}}_i$ such that
$$D_{j_i} \cap CB_i \neq \emptyset \mbox{ and } {\rm{diam}}(D_{j_i}) > {\rm{diam}}(CB_i).$$ We notice
that
$$5CB_i \subset 14 D_{j_i} \mbox{ and } {\rm{diam}}(5CB_i) < 14 \, {\rm{diam}}(D_{j_i}).$$ Therefore
\begin{eqnarray} \label{upper bound for Lambdas of B}
\Lambda_r^{\infty}(\bigcup_{i \in {\mathcal{B}}} 5CB_i) &\le& \sum_{j \in
{\mathcal{J}}} \mbox{diam}(14D_j)^r \\
&\le& 28^r \sum_{j \in {\mathcal{J}}} \mbox{diam}(D_j)^r < 2^{r+1} \, 28^r
\Lambda_r^{\infty}(E). \nonumber
\end{eqnarray}
We let
$${\mathcal{G}}_1 = \{ i \in {\mathcal{G}}: \sum_{j \in {\mathcal{F}}_i} {\rm{diam}}(D_j)^r
<c_0 \, {\rm{diam}}(CB_i)^r  \} \mbox{ and } {\mathcal{G}}_2={\mathcal{G}} \setminus {\mathcal{G}}_1$$
for some $c_0$ to be chosen later. We want to evaluate
$$\sum_{i \in {\mathcal{G}}_2 }  \mbox{diam}(5CB_i)^r.$$
Before we do that, we notice that there exists a number $M$
depending only on the data of $X$ such that every ball $D_j$
intersects at most $M$ pairwise disjoint balls $CB_i$ of bigger
diameter. Therefore
\begin{eqnarray} \label{upper bound for Lambdas of G2} \hspace{8mm}
\Lambda_r^{\infty}(\bigcup_{i \in {\mathcal{G}}_2} 5CB_i) &\le& \sum_{i
\in {\mathcal{G}}_2} \mbox{diam}(5CB_i)^r \le 10^{r} \sum_{i \in
{\mathcal{G}}_2} \mbox{diam}(CB_i)^r\\
&\le& c_0^{-1} \, 10^{r} \sum_{i \in {\mathcal{G}}_2} \left(\sum_{j \in
{\mathcal{F}}_i} \mbox{diam}(D_j)^r\right) \le c_0^{-1} \, 10^{r} \sum_{i
\in {\mathcal{G}}} \left( \sum_{j \in {\mathcal{F}}_i} \mbox{diam}(D_j)^r \right) \nonumber \\
&\le& c_0^{-1} \, M \, 10^{r} \sum_{j \in {\mathcal{J}}}
\mbox{diam}(D_j)^r < c_0^{-1} \, M \, 2^{r+1} 10^r \Lambda_r^{\infty}(E)
\nonumber.
\end{eqnarray}
We show now that if $c_0$ is taken small enough, then
$$\mu(E \cap B_i) \le \frac{1}{2} \mu(B_i) \mbox{ for every } i \in {\mathcal{G}}_1.$$
Indeed, for all $i \in {\mathcal{G}}_1$ we have
\begin{eqnarray*}
\mu(E \cap B_i) &\le& \mu(\bigcup_{j \in {\mathcal{F}}_i} D_j \cap
CB_i) \le \sum_{j \in {\mathcal{F}}_i} \mu(D_j) \le c_{\mu} \sum_{j
\in {\mathcal{F}}_i} {\rm{diam}}(D_j)^{Q}\\ &\le& c_{\mu}
{\rm{diam}}(CB_i)^{Q-r} \sum_{j \in {\mathcal{F}}_i}
{\rm{diam}}(D_j)^{r} \le c_0 c_{\mu} {\rm{diam}}(CB_i)^{Q}.
\end{eqnarray*}
So, if we let $c_0=\frac{1}{2} c_{\mu}^{-2} (2C)^{-Q},$ we get
$$\mu(E \cap B_i) \le \frac{1}{2} \mu(B_i) \mbox{ for every } i \in {\mathcal{G}}_1.$$
From (\ref{upper bound for Lambdas of I}),
(\ref{upper bound for Lambdas of B}), (\ref{upper bound for Lambdas
of G2}), the subadditivity of $\Lambda_r^{\infty},$ and the fact
that ${\mathcal{I}}={\mathcal{G}}_1 \cup {\mathcal{G}}_2
\cup{\mathcal{B}},$ it follows that
$$d(x,y)^r \le \Lambda_s^{\infty}(\bigcup_{i \in {\mathcal{I}}} 5CB_{i}) <
\sum_{i \in {\mathcal{G}}_1}  \mbox{diam}(5CB_i)^r +  2^{r+1} \,
(28^r+c_0^{-1} \,M \, 10^r) \Lambda_r^{\infty}(E).$$ If we choose $c_1$
such that $2^{r+1}(28^r+c_0^{-1} \,M \, 10^r)\,c_1=1-2^{-r},$ then we
notice that
$$(10C)^r \, \sum_{i \in {\mathcal{G}}_1}  \mbox{diam}(B_i)^r \ge
\sum_{i \in {\mathcal{G}}_1}  \mbox{diam}(5CB_i)^r > 2^{-r}
d(x,y)^r$$ whenever $\Lambda_r^{\infty}(E) < c_1 d(x,y)^r.$ Since
${\mathcal{G}}_1 \subset {\mathcal{G}}_0,$ this finishes the proof.
\end{pf}

\begin{Lemma} \label{upper bdd for rel Sobolev cap by fns who are one on the set}
Suppose $1<s,Q<\infty.$ Suppose that $(X,d,\mu)$ is a complete and unbounded Ahlfors $Q$-regular metric measure space
that satisfies a weak $(1,s)$-Poincar\'{e} inequality. Let $\Omega \subset \subset X$ be open and let $E \subset \Omega.$ Suppose
$u \in N_{0}^{1,s}(\Omega)$ is compactly supported in $\Omega.$ If $u \ge 1$ on $E,$ then
$${\mathrm{cap}}_s(E, \Omega) \le \int_{\Omega} g_{u}^s \, d\mu.$$

\end{Lemma}

\begin{pf} Since $u \in N_{0}^{1,s}(\Omega)$ is compactly supported in $\Omega,$ there exists a sequence $\varphi_j \in Lip_{0}(\Omega)$
converging to $u$ in $N^{1,s}(X).$ Without loss of generality we can assume that all the functions $\varphi_j$ are supported in an open set
$U \subset \subset \Omega$ and that the sequence $\varphi_j$ converges to $u$ pointwise $\mu$-a.e. Since $\varphi_j$ is a Cauchy sequence
in $N_{0}^{1,s}(\Omega),$ there is a subsequence, denoted again by $\varphi_j,$ such that
$$||\varphi_j -\varphi_{j+1}||_{N^{1,s}(X)}<2^{-2j} \mbox{ for every $j \ge 1$}.$$
For the open set $$E_{j}=\{ x \in X: |\varphi_j(x)-\varphi_{j+1}(x)| > 2^{-j} \}$$
we have
$$C_{s}(E_{j}) \le 2^{js} ||\varphi_j-\varphi_{j+1}||_{N^{1,s}(X)}^s < 2^{-js}.$$
If we put $$G_j=\bigcup_{k=j} E_{k},$$
we have from the subadditivity of the $s$-capacity that
$$C_s(G_j)^{1/s} \le \sum_{k=j} C_s(E_k)^{1/s} \le \sum_{k=j}^{\infty} 2^{-k} = 2^{1-j}.$$
Thus the sequence $\varphi_j$ converges uniformly outside open sets of arbitrarily small $s$-capacity to a quasicontinuous function $v$ and
we can assume without loss of generality that $v=0$ on $X \setminus U.$
Moreover, $v \in N_{0}^{1,s}(\Omega)$ because $N_{0}^{1,s}(\Omega)$ is a Banach space. On the other hand, $\varphi_j$ converges
to $u$ $\mu$-a.e.\ in $X.$ Thus $u$ and $v$ are two functions in $N^{1,s}(X)$ that agree $\mu$-a.e., hence they agree
$s$-q.e. on $X.$ We let
$$E_0=\{ x \in X: u(x) \neq v(x) \} \mbox{ and } E_1=E \setminus E_0.$$
We fix $\varepsilon \in (0,1).$ We choose an open set $G \subset U$ such that $C_s(G)<\varepsilon$ and $\varphi_j \rightarrow v$
uniformly on $X \setminus G.$
We let $$\widetilde{G}_j = \{ x \in X: \varphi_j(x) > 1-\varepsilon  \}.$$
Then $\widetilde{G}_j$ is open and
$$E_1 \setminus G \subset \widetilde{G}_j \mbox{ for $j \ge j_{\varepsilon}$}.$$
Consequently, for $j \ge j_{\varepsilon}$ we have via the subadditivity of the relative $s$-capacity (see \cite[Theorem 3.2 (vi)]{Cos3})
$${\mathrm{cap}}_{s}(E, \Omega)={\mathrm{cap}}_{s}(E_1, \Omega) \le {\mathrm{cap}}_{s}(\widetilde{G}_j, \Omega) + {\mathrm{cap}}_{s}(G, \Omega).$$
Since $\varphi_j>1-\varepsilon$ on $\widetilde{G}_j,$ we have
$${\mathrm{cap}}_{s}(\widetilde{G}_j, \Omega) \le (1-\varepsilon)^{-s} \int_{\Omega} g_{\varphi_j}^s \, d\mu,$$
and hence by letting $j \rightarrow \infty,$ we obtain
$$ {\mathrm{cap}}_{s}(E, \Omega) \le (1-\varepsilon)^{-s} \int_{\Omega} g_{u}^s \, d\mu + \varepsilon.$$
The lemma follows by letting $\varepsilon \rightarrow 0.$

\end{pf}

We prove Theorem \ref{SA-infty weight s<Q} now.

 \begin{pf}  Since $(X,d,\mu)$ satisfies a weak $(1,s)$-Poincar\'{e} inequality, it follows from \cite{HaK} that
 $(X, d, \mu)$ satisfies in fact a weak $(s,s)$-Poincar\'{e} inequality with possibly another constant $\widetilde{\lambda}.$

 We have that $u \in N^{1,s}_{loc}(X)$ has a minimal $s$-weak upper gradient $g_{u} \in {{\mathcal L}^{s, Q-s}}(X, \mu),$ hence
 we can assume without loss of generality that $u$ is a Borel $s$-quasicontinuous function. Since $g_{u}$ has small
 ${{\mathcal L}^{s, Q-s}}(X, \mu)$ norm, it follows from (\ref{BMO seminorm of u dominated by Morrey norm of grad u}) that
 $u$ has small $\mathrm{BMO}$-seminorm. Therefore,
 from John-Nirenberg lemma, it follows that $w(x)=e^{Q u(x)}$ is an $A_{\infty}$-density with respect to $\mu$
 for some doubling measure $\nu$ with data depending on $X.$
 That is, (see \cite[Theorem 1.4]{MP}, \cite[Theorem A]{MMNO}, and \cite[Theorem 2.2]{Buc}), there exists
 a constant $C$ depending on $s$ and on data of $X$ such that
 \begin{equation}\label{reverse Jensen and doubling measure}
 \frac{1}{\mu(B)} \int_{B} e^{Q (u(x)-u_{B})} \, d\mu(x) <C \mbox{ and } \int_{2B} w(x) \, d\mu(x) \le
 C \int_{B} w(x) \, d\mu(x)
 \end{equation}
 for every ball $B \subset X.$ We write $d \nu(x)=w(x) \, d\mu(x).$ We recall the
 definition of $\delta_{\nu}$ from (\ref{def quasidistance}).
 We shall show that there exists a constant $C \in (0,1]$ such that
 \begin{equation}\label{distance dominates quasidistance}
 d_{\nu}(x_1,x_2):=\inf \sum_{i=1}^k \nu(B_i)^{1/Q}
 \ge C \delta_{\nu}(x_1,x_2)
 \end{equation}
 for all $x_1,x_2 \in X$, where the infimum is taken over
 finite chains of open balls connecting $x_1$ and $x_2$ satisfying
 \begin{equation}\label{chain def}
 x_1 \in B_1, x_2 \in B_k \mbox { and } B_i \cap B_{i+1} \neq \emptyset
 \mbox { for all } i=1, \ldots ,k-1.
 \end{equation}

 Indeed, (\ref{distance dominates quasidistance}) implies both that $d_{\nu}$ is a
 distance and that is comparable to $\delta_{\nu}$ as required
 in (\ref{quasidistance comparable to a distance}). Towards this end, fix $x_1, x_2 \in X,$
 $x_1\neq x_2.$ Let $\gamma$ be a geodesic segment connecting $x_1$ and $x_2,$ and let $a$ be the midpoint of $\gamma.$
 We denote $R=d(x_1, x_2)$ and $B=B(a,R).$

 Let $\eta \in Lip_{0}(6B)$ be a nonnegative $1/R$-Lipschitz function such that $\eta=1$ on $3B.$ Since $u$ is $s$-quasicontinuous
 and Borel, it follows that $v(x)=\eta(x) \, |u(x)-u_{3B}|$ is a Borel $s$-quasicontinuous function in $N_{0}^{1,s}(6B)$
 compactly supported in $6B.$

 Let $E=\{ x \in 3B: |u(x)-u_{3B}|>1  \}.$ We have that $E$ is a Borel set since $u$ is a Borel function.
 Since $v$ is an $s$-quasicontinuous function in $N_{0}^{1,s}(6B)$ compactly supported in $6B,$ we have from
 Lemma \ref{upper bdd for rel Sobolev cap by fns who are one on the set} that
 \begin{eqnarray*}
 {{\mathrm{cap}}}_{s}(E, 6B) & \le & \int_{6B} g_{v}^s \, d\mu
 \le \int_{6\lambda B} (\eta \, g_{u} + |u-u_{3B}| \, g_{\eta})^s \, d\mu\\
 & \le & C \, \int_{6\widetilde{\lambda}B} g_{u}^s \, d\mu \le C \, (6R)^{Q-s} ||g_{u}||_{{{\mathcal L}^{s,Q-s}}(X, \mu)}^s. \nonumber
 \end{eqnarray*}
 This implies that
 \begin{equation}
 \frac{{{\mathrm{cap}}}_{s}(E, 6B)}{(6R)^{Q-s}}
 \le C \, ||g_{u}||_{{{\mathcal L}^{s,Q-s}}(X, \mu)}^s,
 \end{equation}
 which together with \cite[Theorem 4.4]{Cos3} yields
 \begin{equation}\label{rel h, rel cap , grad u}
 \frac{\Lambda_{1}^{\infty}(E)}{R}
 \le C \, \frac{{{\mathrm{cap}}}_{s}(E, 6B)}{(6R)^{Q-s}}
 \le C_{0} ||g_{u}||_{{{\mathcal L}^{s,Q-s}}(X, \mu)}^s.
 \end{equation}
 We choose $\varepsilon>0$ such that
 $C_{0} \, \varepsilon^s < c_1$ where $c_1$ is the
 constant from (\ref{Proj Lemma Lambdar le c1 distr}) and $C_{0}$ is the constant from the last inequality in (\ref{rel h, rel cap , grad u}).

 Now let $B_1, \ldots, B_k$ be an arbitrary chain of balls connecting $x_1$ and $x_2$
 as in (\ref{chain def}). We assume first that $B_i \subset 3B$ for all $i=1, \ldots, k.$ Let
 ${\mathcal{G}}_0$ be defined like in (\ref{Proj Lemma defn of G0}).
 We have
 \begin{eqnarray}\label{SA-infty condition}
 \hspace{5mm} \sum_{i=1}^k \nu(B_i)^{1/Q} & \ge & \sum_{i \in {\mathcal{G}}_0}
 \nu(B_i)^{1/Q} \ge \sum_{i \in {\mathcal{G}}_0} \nu(B_i \setminus E)^{1/Q} =
 \sum_{i \in {\mathcal{G}}_0} \left(\int_{B_i \setminus E} e^{Q u(x)} \, d\mu(x) \right)^{1/Q}\\
 &\ge& \sum_{i \in {\mathcal{G}}_0} \left(\int_{B_i \setminus E} e^{Q (u_{3B}-1)} \, d\mu(x)\right)^{1/Q}
 = e^{u_{3B}-1} \left(\sum_{i \in {\mathcal{G}}_0} \mu(B_i \setminus E)^{1/Q} \right) \nonumber\\
 &\ge& e^{u_{3B}-1} \sum_{i \in {\mathcal{G}}_0} \left(\frac{1}{2} \mu(B_i)\right)^{1/Q} \ge C e^{u_{3B}} \mu(3B)^{1/Q}. \nonumber
 \end{eqnarray}
 From (\ref{reverse Jensen and doubling measure}), (\ref{SA-infty condition}), and the definition of $\delta_{\nu}$
 there exists $C$ such that
 \begin{equation*}
 \sum_{i=1}^k \nu(B_i)^{1/Q} \ge C \left(\int_{3B} e^{Q u(x)} \, d\mu(x)
 \right)^{1/Q} \ge C \delta_{\nu}(x_1,x_2).
 \end{equation*}

 Next, if the chain $(B_i)$ does not lie entirely in $3B,$ then there exists
 a smallest number $k'$ with $1 \le k' \le k$ such that
 $B_{k'} \cap S(a, 2R) \neq \emptyset.$ Let $x_0 \in B_{k'} \cap S(a, 2R).$ Then $B_1, \ldots, B_{k'}$
 is a chain of balls connecting $x_1$ and $x_0$ and $d(x_1,x_2) \le d(x_1,x_0).$

 If $B_{k'} \subset 3B,$ then from the fact that $x_1 \in B_1 \cap 2B$ and from the definition of $k'$
 it follows that the subchain $B_1, \ldots, B_{k'}$ is contained in $3B.$ Therefore we can apply the
 preceding argument to the chain $B_1, \ldots, B_{k'}$ connecting the points $x_1$ and $x_0$ to conclude that
 (\ref{distance dominates quasidistance}) holds; in the opposite case, $\mbox{ diam } B_{k'} \ge R.$ The doubling
 condition for $\nu$ then implies $\nu(B) \le C \nu(B_{k'}).$
 Thus, (\ref{distance dominates quasidistance}) is true in all
 cases. This finishes the proof.
 \end{pf}

 \begin{Lemma} \label{upper bdd for rel Besov cap by fns who are one on the set}
 Suppose $1<Q<p<\infty.$ Suppose $(X,d,\mu)$ is a complete and unbounded Ahlfors $Q$-regular metric measure space.
 Let $\Omega \subset \subset X$ be open and let $E \subset \Omega.$ Suppose $u \in B_p^{0}(\Omega)$ is a $B_p$-quasicontinuous function
 compactly supported in $\Omega.$ If $u \ge 1$ on $E,$ then
 $${\mathrm{cap}}_{B_p}(E, \Omega) \le [u]_{B_p(\Omega)}^p.$$

 \end{Lemma}

 \begin{pf} Since $u \in B_p^{0}(\Omega)$ is compactly supported in $\Omega,$ there exists a sequence $\varphi_j \in Lip_{0}(\Omega)$
converging to $u$ in $B_p(X).$ (See \cite[Lemma 3.14]{Cos2}.) Without loss of generality we can assume that all the functions $\varphi_j$
are supported in an open set $U \subset \subset \Omega$ and that the sequence $\varphi_j$ converges to $u$ pointwise $\mu$-a.e.
Since $\varphi_j$ is a Cauchy sequence
in $B_p^{0}(\Omega),$ there is a subsequence, denoted again by $\varphi_j,$ such that
$$||\varphi_j -\varphi_{j+1}||_{L^{p}(X)}^p+[\varphi_j -\varphi_{j+1}]_{B_p(X)}^p<2^{-j(p+1)} \mbox{ for every $j \ge 1$}.$$
For the open set $$E_{j}=\{ x \in X: |\varphi_j(x)-\varphi_{j+1}(x)| > 2^{-j} \}$$
we have
$${\mathrm{Cap}}_{B_p}(E_{j}) \le 2^{jp} \left( ||\varphi_j -\varphi_{j+1}||_{L^{p}(X)}^p+[\varphi_j -\varphi_{j+1}]_{B_p(X)}^p \right)  < 2^{-j}.$$
If we put $$G_j=\bigcup_{k=j} E_{k},$$
we have via the subadditivity of the Besov $p$-capacity (see \cite[Theorem 5.3 (ii)]{Cos2})
$${\mathrm{Cap}}_{B_p}(G_j) \le \sum_{k=j} {\mathrm{Cap}}_{B_p}(E_k) \le \sum_{k=j}^{\infty} 2^{-k} = 2^{1-j}.$$
Thus the sequence $\varphi_j$ converges uniformly outside open sets of arbitrarily small Besov $p$-capacity to a $B_p$-quasicontinuous function $v$ and we can assume without loss of generality that $v=0$ on $X \setminus U.$
Moreover, $v \in B_p^{0}(\Omega)$ because $B_p^{0}(\Omega)$ is a Banach space. On the other hand, $\varphi_j$ converges
to $u$ $\mu$-a.e.\ in $X.$ Thus $u$ and $v$ are two $B_p$-quasicontinuous functions in $B_p(X)$ that agree $\mu$-a.e., hence they agree $B_p$-q.e. (See \cite[p.262]{Kil} and \cite[Theorem 5.16]{Cos2}.)
We let
$$E_0=\{ x \in X: u(x) \neq v(x) \} \mbox{ and } E_1=E \setminus E_0.$$
We fix $\varepsilon \in (0,1).$ We choose an open set $G \subset U$ such that ${\mathrm{Cap}}_{B_p}(G)<\varepsilon$ and $\varphi_j \rightarrow v$
uniformly on $X \setminus G.$
We let $$\widetilde{G}_j = \{ x \in X: \varphi_j(x) > 1-\varepsilon  \}.$$
Then $\widetilde{G}_j$ is open and
$$E_1 \setminus G \subset \widetilde{G}_j \mbox{ for $j \ge j_{\varepsilon}$}.$$
Consequently, for $j \ge j_{\varepsilon}$ we obtain via the subadditivity of the relative Besov $p$-capacity (see \cite[Theorem 4.2 (vi)]{Cos2})
$${\mathrm{cap}}_{B_p}(E, \Omega)={\mathrm{cap}}_{B_p}(E_1, \Omega) \le {\mathrm{cap}}_{B_p}(\widetilde{G}_j, \Omega) + {\mathrm{cap}}_{B_p}(G, \Omega).$$
Since $\varphi_j>1-\varepsilon$ on $\widetilde{G}_j,$ we have
$${\mathrm{cap}}_{B_p}(\widetilde{G}_j, \Omega) \le (1-\varepsilon)^{-p} [\varphi_j]_{B_p(\Omega)}^p,$$
and hence by letting $j \rightarrow \infty,$ we obtain
$$ {\mathrm{cap}}_{B_p}(E, \Omega) \le (1-\varepsilon)^{-p} [u]_{B_p(\Omega)}^p + \varepsilon.$$
The lemma follows by letting $\varepsilon \rightarrow 0.$

 \end{pf}

 Now we prove Theorem \ref{SA-infty weight p>Q}.
 \begin{pf} Since $u$ has small Besov $p$-seminorm, it follows via
 (\ref{BMO seminorm of u dominated by Besov seminorm of u}) that
 $u$ has small $\mathrm{BMO}$-seminorm. Therefore, by John-Nirenberg lemma
 (see \cite[Theorem 1.4]{MP}, \cite[Theorem A]{MMNO}, and \cite[Theorem 2.2]{Buc}) there exists a
 constant $C$ depending on $p$ and on the data of $X$ such that $w=e^{Qu}$ is an $A_{\infty}$-density with respect to $\mu$
 for some doubling measure $\nu$ satisfying (\ref{reverse Jensen and doubling measure}) with $C.$
 We write $d \nu(x)=w(x) \, d\mu(x).$ We recall the definition of $\delta_{\nu}$ from
 (\ref{def quasidistance}).
 We shall show that there exists a constant $C \in (0,1]$
 such that (\ref{distance dominates quasidistance}) holds
 for all $x_1,x_2 \in X,$ where the infimum is taken over
 finite chains of open balls connecting $x_1$ and $x_2$ satisfying (\ref{chain def}).

 Indeed, (\ref{distance dominates quasidistance}) implies both that $d_{\nu}$ is a
 distance and that is comparable to $\delta_{\nu}$ as required
 in (\ref{quasidistance comparable to a distance}). Towards this end, fix $x_1, x_2 \in X.$
 We can assume without loss of generality that $x_1 \neq x_2.$ Let $\gamma$ be a geodesic segment connecting $x_1$ and $x_2,$ and
 let $a$ be the midpoint of $\gamma.$ We denote $R=d(x_1, x_2)$ and $B=B(a,R).$

 Let $\eta \in Lip_0(6B)$ be a nonnegative $1/R$-Lipschitz function such that $\eta=1$ on $3B.$ We let the function $v$ be
 defined by $v=\eta |u-u_{3B}|.$ Since $u \in L_{loc}^{1}(X)$ with $[u]_{B_p(X)}<\infty,$ it follows via \cite[Lemma 3.11]{Cos2} that
 $v \in B_p^{0}(6B).$ Moreover, since $v$ is compactly supported in $\Omega,$ it follows via \cite[Theorem 5.12]{Cos2} that there exists
 a $B_p$-quasicontinouus function $\widetilde v$ compactly supported in $\Omega$ such that $v=\widetilde v$ $\mu$-a.e. in $\Omega.$
 Moreover, we can assume that $\widetilde v$ is Borel via an argument similar to the one from Lemma
 \ref{upper bdd for rel Besov cap by fns who are one on the set}. Thus we can assume without loss of generality that $v$ is
 $B_p$-quasicontinuous and Borel.

 Let $$E=\{ x \in 3B: |u(x)-u_{3B}|>1  \}=\{ x \in 3B: v(x)>1 \}.$$ Hence $E$ is a
 Borel set since $v$ is a Borel function. From Lemma \ref{upper bdd for rel Besov cap by fns who are one on the set} and
 \cite[Lemma 3.11]{Cos2} we obtain
 \begin{eqnarray*}
 {\mathrm{cap}}_{B_p}(E, 6B) &\le& [v]_{B_p(6B)}^p \le C \, [u]_{B_p(X)}^p
 \end{eqnarray*}
 where $C$ is a constant that depends only on $p$ and on data of $X.$
 This together with \cite[Theorem 4.10]{Cos2} yields
 \begin{equation}\label{rel h, rel cap , u B_p}
 \frac{\Lambda_{1}^{\infty}(E)}{R} \le C \, {{\mathrm{cap}}}_{B_p}(E, 6B) \le C \, [v]_{B_p(6B)}^p \le C_0 \, [u]_{B_p(X)}^p,
 \end{equation}
 where $R$ is the radius of $B$ and $C_{0}$ is a constant that depends only on $p$ and on data of $X.$
 We choose $\varepsilon>0$ such that $C_{0} \, \varepsilon^p < c_1$ where $c_1$ is the
 constant from (\ref{Proj Lemma Lambdar le c1 distr}) and $C_{0}$ is the constant
 from the last inequality in (\ref{rel h, rel cap , u B_p}).

 Now let $B_1, \ldots, B_k$ be an arbitrary chain of balls connecting $x_1$ and $x_2$
 as in (\ref{chain def}). We assume first that $B_i \subset 3B$ for all $i=1, \ldots, k.$
 Let ${\mathcal{G}}_0$ be defined like in (\ref{Proj Lemma defn of G0}).
 The proof now continues like in Theorem \ref{SA-infty weight s<Q}, with
 the only difference that the constants who depended on $s$ and on the data of $X$ will now depend
 on $p$ and on the data of $X.$

 \end{pf}

 Theorem \ref{SA-infty weight s<Q} yields the following consequence:

 \begin{Theorem}\label{Resh equivalence for s<Q} Let $1<s \le Q<\infty$ be fixed. We assume that $s>Q-1.$
 Suppose $(X,d,\mu)$ is a metric measure space as in Theorem \ref{SA-infty weight s<Q}. Let $u$ be a Borel function in
 $N_{loc}^{1,s}(X)$ such that it has a minimal $s$-weak upper gradient $g_{u}$ in the Morrey space ${{\mathcal L}^{s, Q-s}}(X, \mu).$
 There exists a constant $\varepsilon>0$ depending only on $s$ and on the data of $X$ such that if
 $$||g_{u}||_{{{\mathcal L}^{s, Q-s}}(X, \mu)} < \varepsilon,$$ then
 \begin{equation} \label{Resh dominance} \delta_{\nu}(x_1,x_2) \le C D_{\nu}(x_1,x_2) \mbox { for all  $x_1,x_2$ in } X,
 \end{equation}
 where $C>0$ is a constant depending only on $s$ and on the data associated with $X$ and
 \begin{equation} \label{definition of Dnu}
 D_{\nu}(x,y)=\inf \bigg \lbrace \int_{\gamma} e^{u} ds: \mbox{ $\gamma$ a rectifiable curve connecting $x,y$} \bigg \rbrace.
 \end{equation}
 \end{Theorem}

 For a discussion about line integration see \cite[Chapter 7]{Hei}.

 Theorem \ref{SA-infty weight p>Q} yields the following consequence:

 \begin{Theorem}\label{Resh equivalence for p>Q} Let $1<s<Q<p<\infty$ be fixed.
 Suppose $(X,d,\mu)$ is a metric measure space as in Theorem \ref{SA-infty weight p>Q}.
 Let $u$ be in $L^{1}_{loc}(X)$ such that $[u]_{B_p(X)}<\infty.$
 There exists a constant $\varepsilon>0$ depending only on $p$ and on the data of $X$ such that if
 $[u]_{B_p(X)} < \varepsilon,$ then
 \begin{equation*} \delta_{\nu}(x_1,x_2) \le C D_{\nu}(x_1,x_2) \mbox { for all  $x_1,x_2$ in } X,
 \end{equation*}
 where $C>0$ is a constant depending only on $p$ and on the data associated with $X$ and $D_{\nu}$ is defined as in (\ref{definition of Dnu}).
 Here $\widetilde u$ is a $B_p$-quasicontinuous Borel representative of $u.$
 \end{Theorem}

 One should compare the metrics $D_{\mu}$ in Theorems \ref{Resh equivalence for s<Q}
 and \ref{Resh equivalence for p>Q} to those studied in \cite{BL}, \cite{Res} and \cite{Cos1}.
 \begin{Question} Another open question is whether or not the inequality (\ref{Resh dominance})
 can be reversed in general. The answer is yes in the Euclidean setting when $n \ge 2.$ (See \cite[Theorems 5.4 and 5.5]{Cos1}.)
 \end{Question}

 We prove now Theorem \ref{Resh equivalence for s<Q}.
 \begin{pf}  It is easy to see that $D_{\nu}$ is indeed symmetric, nonnegative and
 satisfies the triangle inequality. From (\ref{Resh dominance}) it would
 follow immediately that $D_{\nu}$ is a distance function dominating
 $\delta_{\nu}.$ So fix $x_1,x_2$ in $X.$ We can assume without loss
 of generality that $x_1\neq x_2.$ Like before, let $a$ be a point such that $d(x_1,a)=d(a,x_2)=R/2,$
 where $R=d(x_1,x_2).$ We denote $B=B(a,R).$
 Like in the proof of Theorem \ref{SA-infty weight s<Q}, let
 $$v=\eta |u-u_{3B}| \mbox{ and } E=\{ x \in 3B: |u(x)-u_{3B}|>1  \},$$ where $\eta \in Lip_{0}(6B)$ is a
 nonnegative $1/R$-Lipschitz function such that $\eta=1$ on $3B.$ We notice that $E$ is a Borel set
 and $v$ is a Borel and $s$-quasicontinuous function in $N_{0}^{1,s}(6B)$ compactly supported
 in $6B.$

 Let $\gamma$ be a rectifiable curve connecting $x_1$ and $x_2$ and let $|\gamma|$ be its image. We
 assume first that $|\gamma| \subset 3B.$ We obviously have
 \begin{equation}\label{ineq rect curve}
 \int_{\gamma} e^u ds \ge \int_{\gamma \cap (3B \setminus E)} e^u ds.
 \end{equation}

 As in the proof of Theorem \ref{SA-infty weight s<Q}, we have
 \begin{eqnarray*}
 \frac{\Lambda_{1}^{\infty}(E)}{R}
 \le C \, \frac{{{\mathrm{cap}}}_{s}(E, 6B)}{(6R)^{Q-s}}
 \le C_{0} ||g_u||_{{\mathcal L}^{s,Q-s}(X, \mu)}^s<C_{0} \, \varepsilon^s,
 \end{eqnarray*}
 hence
 \begin{eqnarray*}
 \Lambda_{1}^{\infty}(|\gamma| \cap (3B \setminus E)) &\ge& \Lambda_{1}^{\infty}(|\gamma| \cap 3B)-\Lambda_{1}^{\infty}(|\gamma| \cap E) \\
 &\ge& R- \Lambda_{1}^{\infty}(E)\ge (1-c_1)R
 \end{eqnarray*}
 if $\varepsilon>0$ is small enough, where $c_1$ is the constant from (\ref{Proj Lemma Lambdar le c1 distr}).

 Thus we obtain
 \begin{eqnarray}\label{ineq rect curve 1}
 \int_{\gamma} e^u ds & \ge & \int_{\gamma \cap (3B \setminus E)} e^u ds
 \ge \int_{\gamma \cap (3B \setminus E)} e^{u_{3B}-1} ds\\
 &\ge& \Lambda_{1}^{\infty}(|\gamma| \cap (3B \setminus E)) e^{u_{3B}-1}
 \ge C \, R \, e^{u_{3B}}\nonumber\\
 & \ge & C \left(\int_{3B} e^{Q u_{3B}} d\mu \right)^{1/Q} \ge
 C \left(\int_{3B} e^{Q u(z)} d\mu(z) \right)^{1/Q},  \nonumber
 \end{eqnarray}
 where the last inequality follows from (\ref{reverse Jensen and doubling measure}).
 Therefore
 \begin{equation}\label{ineq rect curve 2}
 \int_{\gamma} e^u ds \ge C \left(\int_{3B} e^{Q u(z)} \, d\mu(z)
 \right)^{1/Q} \ge C \delta_{\nu}(x_1,x_2) \mbox{ whenever } |\gamma| \subset 3B.
 \end{equation}

 Now we assume that $|\gamma| \setminus 3B \neq \emptyset.$ Suppose that $\gamma$ is parametrized
 by its arc length parametrization. Let $t_0=\inf \{ t \in [0, l_{\gamma}]: \gamma(t) \notin B(a,2R) \}.$
 Then, since $\gamma$ is a path with $\gamma(0), \, \gamma(l_{\gamma}) \in 2B,$ it follows that
 $$0<t_0<l_{\gamma} \mbox{ and } \gamma([0, t_0]) \subset \overline{B}(a,2R).$$
 Let $x_0=\gamma(t_0)$ and let $\widetilde \gamma$ be the restriction of $\gamma$ to $[0, t_0].$ Then $x_0 \in S(a, 2R)$
 and $d(x_1, x_2) \le d(x_1, x_0).$ Therefore
 \begin{eqnarray*}
 \Lambda_{1}^{\infty}(|\gamma| \cap (3B \setminus E)) &\ge& \Lambda_{1}^{\infty}(|\widetilde \gamma| \cap (3B \setminus E)) \ge \Lambda_{1}^{\infty}(|\widetilde \gamma|) -\Lambda_{1}^{\infty}(|\widetilde \gamma| \cap E) \\
 &\ge& d(x_1, x_0) - \Lambda_{1}^{\infty}(E) \ge (1-c_1)R
 \end{eqnarray*}
 if $\varepsilon>0$ is small enough, where $c_1$ is the constant from (\ref{Proj Lemma Lambdar le c1 distr}).
 By repeating the argument from (\ref{ineq rect curve 1}) with $\widetilde \gamma$ instead of $\gamma,$ we obtain
 \begin{equation*}
 \int_{\widetilde \gamma} e^{u} ds \ge C \left(\int_{3B} e^{Q u(z)} \, d\mu(z) \right)^{1/Q}.
 \end{equation*}
 The desired conclusion follows.

 \end{pf}

  Now we prove Theorem \ref{Resh equivalence for p>Q}.
 \begin{pf} For the existence and "uniqueness" of $B_p$-quasicontinuous Borel representatives of $u$ see \cite[Corollary 5.19]{Cos2}.
 We notice that $D_{\nu}$ does not depend on the choice of the $B_p$-quasicontinuous Borel representative. Indeed,
 if $\widetilde u$ and $\widetilde v$ are two such representatives, then from \cite[Corollary 5.18]{Cos2}
 we have $\widetilde u=\widetilde v$ $B_p$-q.e., which implies via \cite[Corollary 4.15]{Cos2} that
 \begin{equation*}
 \int_{\gamma} e^{\widetilde u} ds=\int_{\gamma} e^{\widetilde v} ds
 \end{equation*}
 for every rectifiable curve $\gamma$ in $X.$ So we can assume without loss of generality
 that $u$ is a $B_p$-quasicontinuous Borel function itself.

 It is easy to see that $D_{\nu}$ is indeed symmetric, nonnegative and
 satisfies the triangle inequality. From (\ref{Resh dominance}) it would
 follow immediately that $D_{\nu}$ is a distance function dominating
 $\delta_{\nu}.$ So fix $x_1,x_2$ in $X.$ We can assume without loss
 of generality that $x_1\neq x_2.$ Like before, let $a$ be a point such that $d(x_1,a)=d(a,x_2)=R/2,$
 where $R=d(x_1,x_2).$ We denote $B=B(a,R).$

 Let $\eta \in Lip_0(6B)$ be a nonnegative $1/R$-Lipschitz function such that $\eta=1$ on $3B.$
 Like in the proof of Theorem \ref{SA-infty weight p>Q}, let
 $$v=\eta |u-u_{3B}| \mbox{ and } E=\{ x \in 3B: |u(x)-u_{3B}|>1  \}= \{ x \in 3B: v(x)>1 \}.$$ We notice that $E$ is a Borel set
 and $v$ is a $B_p$-quasicontinuous Borel function in $B_p^{0}(6B)$ compactly supported
 in $6B.$

 Let $\gamma$ be a rectifiable curve connecting $x_1$ and $x_2$ and let $|\gamma|$ be its image. We
 assume first that $|\gamma| \subset 3B.$ We obviously have
 \begin{equation}\label{ineq rect curve}
 \int_{\gamma} e^u ds \ge \int_{\gamma \cap (3B \setminus E)} e^u ds.
 \end{equation}

 Like in the proof of Theorem \ref{SA-infty weight p>Q}, we have
 \begin{eqnarray*}
 \frac{\Lambda_{1}^{\infty}(E)}{R}
 \le C \, {{\mathrm{cap}}}_{B_p}(E, 6B)
 \le C_{0} [u]_{B_p(X)}^p<C_{0} \, \varepsilon^p.
 \end{eqnarray*}
 The proof now continues like in Theorem \ref{Resh equivalence for s<Q},
 with the only difference that the constants who depended on $s$ and on the data of $X$ will now depend
 on $p$ and on the data of $X.$

 \end{pf}

%\bibliography{Bibliopapers}
\end{document}